\documentclass[11pt]{article}
\usepackage[ansinew]{inputenc}
\usepackage[psamsfonts]{amssymb}
\def\bk{{l\!k}}\title{Configuration space and Massey products}
\author{Yves Felix and Jean-Claude Thomas }

\begin{document}
\maketitle
\abstract{The purpose of this paper is to study and compare the collapsing  of
two
 spectral sequences converging to the cohomology of  a configuration
 space. The non collapsing of these spectral sequences is related,
 in some cases, to the existence of Massey products in the cohomology of the
manifold $M$.}

 \vspace{1cm}
Let $\bk$ be a field and  $M$ be an oriented connected closed $m$-dimensional
real manifold. We consider the configuration space of ordered $n$-tuples of
points in $M$
$$
F(M,n)= \{ (x_1, \ldots, x_n)\in M^{ n} \,\vert\,  x_i \neq
x_j\,, \mbox{if} \, i\neq j\, \}\,.
$$
The symmetric group $\Sigma_n$ acts freely on $F(M,n)$ and the quotient
manifold
$B(M,n) =
F(M,n)/\Sigma_n$ is called the space of unordered $n$-tuples of distinct
points
in $M$.

 The determination of the homology of the configuration spaces
 $F(M,n)$ and $B(M,n)$         has been at the origin of a lot of works:
 In \cite{Bo},         Bodigheimer, Cohen and Taylor       compute         the
 homology groups of $B(M,n)$ when $M$ is an odd-dimensional closed
 manifold. In \cite{LoM}, L\"offler and Milgram       compute         the
 $\mathbb Z/2$-homology groups of $B(M,n)$ for any closed manifold
 $M$. The answers depend only on the homology groups of $M$. In
 \cite{FT}, we have computed the rational homology of $B(M,n)$: the
 answer depends on the cohomology algebra of $M$, and not only on
 the rational Betti numbers when $M$ is even-dimensional.

There is no direct procedure to compute the cohomology algebra
$H^*(F(M,n);\bk)$ in general.
There are however two   spectral sequences converging  to   $H^*(F(M,n);\bk)$.
 The first one, we call the Cohen-Taylor spectral sequence for
$H^*(F(M,n);\bk)$, is the Leray
spectral sequence  for the  inclusion
$F(M,n)
\hookrightarrow M^n$, and converges as an algebra to $H^*(F(M,n);\bk)$
(\cite{C}). Its
$E_1$ term   is the quotient of
the  bigraded
commutative algebra
 $$
 \, H^{\otimes n} \otimes \land (x_{st})
$$
 where

\vspace{1mm}\begin{minipage}{12.5cm}$H^i=H^i(M;\bk)$ is concentrated in
bidegree $(0,i)$ and
 $\land (x_{st}) $ denotes the free bigraded commutative algebra generated
by the
 elements $ x_{st}$ of bidegree $(0,m-1)$ for $ 1\leq s,t\leq n \,, \, s \neq
t$\end{minipage}

\vspace{1mm}\noindent  by the ideal  $I$  generated by the elements
$$
(e_s-e_t)(a) \otimes x_{st} \,,\quad  x_{st}^2\,, \quad     x_{st} -
(-1)^mx_{ts}\,,\quad
x_{st}x_{su} + x_{su}x_{st} + x_{us}x_{ut}\,,
$$
where   $e_i(a) = 1 \otimes
\cdots \otimes a \otimes 1 \cdots \otimes 1$ ($a$ appears in the $i^{th}$
position).

   The differential $d_1$ is given by
$$
d_1(x_{st})= p^\ast _{s,t}(\delta ^M)
$$
where $p_{s,t} ; M^{ n} \to M^{ 2}$ denote the   projection on the
$s^{\scriptstyle th}$ and the $t^{\scriptstyle th}$ factors
and $\delta ^M \in
H^m(M^{2} ;\bk)= \left( H \otimes H\right)^m$ denotes the diagonal class.

 Using a
weight-filtration argument,   Totaro (\cite{T})  has proved  that for   a
smooth
projective variety $M$,
the differential  $d_1$ is the only non trivial differential and  that
$H(E_1,d_1)\cong
H^*(F(M,n);\bk)$ as an algebra (see also \cite{K}). He also asks for an
example of
a manifold
for which the Cohen-Taylor spectral sequence has  other nonzero differentials.

The second spectral sequence, we call the Bendersky-Gitler spectral
sequence for
$H^*(F(M,n);\bk)$, is
obtained by considering  the nerve of the covering of the  fat diagonal
 $$
D^nM = \{ \, (x_1, \ldots , x_n) \in M^{ n}\,\vert\, x_i = x_j \,,
 \mbox{ some pair } i\neq j\, \}
$$
  by  the subspaces
 $$
D_{ij} = \{\, (x_1, \ldots , x_n\,) \in M^{ n}\,\vert \, x_i = x_j\,\}\,.
$$
This defines a double cochain complex $
{\cal C}^\ast(X_\bullet) : {\cal C}^\ast (X_0) \stackrel{d'}{\to}
 {\cal C}^\ast (X_1) \stackrel{d'}{\to} \cdots \to {\cal
 C}^\ast (X_n)$
 with $X_r = \amalg \, ( D_{i_1j_1}\cap \cdots
 \cap D_{i_{r+1}j_{r+1}})\,.$
As proved   by   Bendersky and   Gitler (\cite{BG})  the first   spectral
sequence  produced by this bicomplex    is a spectral sequence of
$H^*(M^n;\bk)$-modules converging to  $H^\ast(M^{ n}, D^nM)$.
By Lefschetz duality we have $H^\ast(M^{ n}, D^nM) \cong H_{mn-\ast}
(F(M,n);\bk
)$. Thus
the dual of the Bendersky-Gitler spectral sequence converges to
$H^{mn-\ast}(F(M,n);\bk )$.

  When $\bk$ is a field of characteristic zero and  $M$ is rationally
formal, $d_1$ is the only non zero
differential. Recall that
a space $X$  is rationally formal if there is a sequence of
 quasi-isomorphisms connecting the de Rham algebra of differential forms on
$M$,
$\Omega^*(M)$,
 and its real
 cohomology, $(H^*(M;\mathbb R),0)$). A smooth projective variety is a formal
space
(\cite{Del}) and if $M$ is a formal space all the Massey products are trivial.
Bendersky and Gitler
have  conjectured that higher differentials are determined by higher order
Massey
products.

Our main results are:

\vspace{3mm}\noindent  {\bf Theorem 1.}

\noindent {\sl The Poincar\'e  duality of $H^*(M;\bk)$
induces
an isomorphism of  $ H^*(M;\bk)^{\otimes n}$-modules between  the $E_2$
term  in    the Cohen-Taylor spectral
sequence   and the dual of the $nm$-suspension of the $E_2$ term in the
Bendersky-Gitler spectral
sequence.}

\vspace{3mm} Since  each $H_{k}(F(M,n);\bk)$ is finite dimensional, Theorem
1 implies that
the
 Bendersky spectral
sequence for $H_{mn-\ast}(F(M,n);\bk)$ collapses at the $E_2$ term if and only
if  the Cohen-Taylor
spectral sequence for $H^*(F(M,n);\bk)$ collapses at the $E_2$ term. In
particular,  the
collapsing theorem of Bendersky-Gitler implies the first part  of the
collapsing
result of Totaro.

\vspace{3mm}\noindent {\bf Theorem 2.} {\sl If $n\leq 3$, then both spectral
 sequences  collapse at the $E_2$ term.}

 \vspace{3mm}\noindent {\bf Theorem 3.} {\sl  Let $\bk$ be   of
characteristic zero and $M$ be   simply
 connected. Suppose there exist indecomposable elements  $a,b,c$
and $d\in
H^*(M;\bk)$, such that
 \begin{enumerate}
 \item[$\bullet$] $ab = ac= ad= bc=bd=cd= 0$,
 \item[$\bullet$]  the triple Massey product $\langle b,c,d\rangle $ is a
 nonzero indecomposable element, and
 \item[$\bullet$] $a$ is not a linear
 combination of the elements $b,c,d$
and $\langle b,c,d\rangle $,
 \end{enumerate}  then
the Cohen-Taylor spectral
 sequence for   $H^*(F(M,4);\bk)$ does not collapse at the $E_2$ term. }

 \vspace{3mm}\noindent {\bf Theorem 4.} {\sl  Let $k$ be   of
characteristic zero, $M$ be   simply
 connected  and $N = M\# ( S^2 \times S^{m-2}) $ be the connected sum of $M$
with
 a product  of two  spheres.  If
  there is  a triple matrix Massey
 product in $H^*(M;\bk)$ represented by an indecomposable
 element, then the Cohen-Taylor spectral
 sequence for   $H^*(F(N,4);\bk)$ does not collapse at the $E_2$ term.  }

\vspace{3mm}
 Theorems 3 and 4    answer  the question of  Totaro quoted above. We give
explicit formulae for the differential $d_2$ of the Bendersky-Gitler
 spectral sequence in section 4.
This enables a lot of   examples. For instance at the
end of section 4 we
construct a nonzero $d_2$ for $F(M,4)$  when   $M$ is  the sphere tangent
bundle
to
the manifold  $S^2\times S^2$.

The  main ingredient in the text     is the  bicomplex
$C(n, A)$
constructed in section 3 for any differential graded commutative algebra $A$
using a convenient
family of oriented graphs with
$n$ vertices  and such that no two edges have the same range. The
associated  spectral sequence   is isomorphic
to the Bendersky-Gitler spectral sequence in characteristic zero.
 The simplicity of the  complex allows explicit computations. In particular
we prove:

 \vspace{3mm}
\noindent{\bf Proposition 5.} {\sl There exists   a short exact sequence
 $$
0 \to (I_M)^* \to H^*(F(M,3);\bk) \to (\mbox{Hom} (\Omega^1_{H/\bk},
\bk))^{*-1+3m} \to 0\,.
$$ Here
 $I_M$ is the image         of $H^*(M^3;\bk )$ along the map induced by
 the injection $F(M,3) \hookrightarrow M^3$ and         $\Omega^1_{H/\bk}$
 is the module of K\"ahler differentials, $H=H^*(M;\bk)$.}

 \vspace{3mm} In the first section of the paper we describe the
 Bendersky-Gitler spectral sequence in terms of graphs. In section 2 we prove
theorem 1.    In section 3 we construct the bicomplex $C(n,A)$ and prove
theorems 2.  Section 4 is
devoted to Massey products, the
 homology of $F(M,4)$ and the proof of Theorems 3 and 4.

 We are  grateful to Martin Markl and Stephan Papadima for having signaled
us a
mistake in the
previous version of this paper.

 \section{The $E_1$  term in the  Bendersky-Gitler spectral sequence}

 To go further in the description of the $E_1$ term of the Bendersky-Gitler
spectral  sequence, we
introduce for any  differential graded commutative algebra $(A,d)$ the
differential  bigraded
 algebra $E(n,A)$ defined as follows: As an an algebra $E(n,A)$ is the
quotient
of
 $$
 A^{\otimes n} \otimes \land (e_{ij})
$$
where

the generators $e_{i,j}\,, \quad 1\leq i<j\leq n $,  have bidegree $(0,1)$

the elements $a_1\otimes ...\otimes a_n \otimes 1$ has bidegree $(\sum
|a_1|,0))$

\noindent by the the ideal   generated by the elements
 $$(e_i-e_j)(a) \otimes e_{ij}$$
 where $e_i(a) = 1 \otimes \cdots \otimes a \otimes 1 \cdots
\otimes 1$, with $a$
 appearing in the $i^{th}$ position.

The differential $d$ is defined by
$$
d = d' + d''
$$
where $d'$ is the multiplication by $ \sum e_{ij}$ and $d''$ is the
differential induced by the differential $d$ on
$A$.

\vspace{3mm}
We are particulary interested  in the  case when  $A$ is a commutative
model of
$M$ over
$\bk$. Such a commutative model exists in the following situations
\begin{enumerate}
\item[$\bullet$]   $M$ is a  $\bk$-formal space, any $\bk$ (with
$A=H^\ast(M;\bk)$),
\cite{El},
\cite{Del},
\item[$\bullet$]    $M$ is $r$-connected and the characteristic $p$ of the
field $\bk$ satisfies $p >
\frac{m}{r}$
(\cite{An}),
\item[$\bullet$]    $\bk= {\mathbb Q}$,  (with $A = A_{PL}(M)$ or any
Sullivan model of
$M$), \cite{Su},
\cite{FHT},
\item[$\bullet$]    $\bk= {\mathbb R}$, (with $A$ = the de Rham  complex of
$M$).
\end{enumerate}

By filtering $E(n,A)$ by the powers of the ideal generated by the $e_{ij}$,
Bendersky and Gitler construct a new spectral sequence.
They prove (cf. Theorems 2.7 and 3.1 in \cite {BG})

\vspace{3mm} \noindent {\bf Theorem A.} {\sl For any field $\bk$,
the cohomology $H^*(E_1,d_1)$ of the Bendersky-Gitler spectral sequence is
isomorphic to $H^*(E(n,H^*(M;\bk)))$.}

\vspace{3mm} \noindent {\bf Theorem B.} {\sl Let  $A$ be  a  commutative model
of $M$ over
$\bk$. The      spectral sequence  associated to
 the double complexe $E(n,A)$ is isomorphic to   the
Bendersky-Gitler spectral
sequence for $H_{mn-\ast}(F(M,n);\bk)$ from the $E_2$ term.}

 \vspace{3mm} Following   Bendersky and Gitler \cite{BG}, it is convenient
 to describe $E(n,A)$ in terms of oriented graphs.

\vspace{3mm}
 Let $n$ be a
 fixed integer and let ${\cal G}(n)$ be the set of graphs $G$ with
 set of vertices $V(G) = \{ 1, 2, \ldots,n\}$ and set of edges $E(G) \subset\{
(i,j)\,, i,j
\in V(G), i<j\}$.
 We note $e_G = e_{i_1j_1}\ldots e_{i_kj_k}$, if $E(G) =\{
 {(i_1,j_1)}\ldots         {(i_k,j_k)}\}$, where the $(i_k,j_k)$ are
ordered by
 lexicographic order.

 We denote by $S_1$, $S_2$, $\ldots $, $S_{l(G)}$ the components of
 $G$ ordered by the smallest vertex, and if  $S_i$ is the ordered set $ \{
i_1,
i_2, \ldots ,
i_r\}$  by $a_{S_i}$ the product $a_{i_1}a_{i_2} \ldots
 a_{i_r} \in A$.

 Then the natural map $a_1\otimes \ldots \otimes
 a_n \otimes e_G \mapsto \varepsilon_G a_{S_1}\otimes \ldots \otimes
 a_{S_{l(G)}} \otimes e_G$ ($\varepsilon_G$ is the  graded signature)
induces an
isomorphism of
graded vector spaces
 $$
\Phi : E(n,A) \to       \oplus_{G\in {\cal G}(n)}
 A^{\otimes l(G)}e_G\, :={\cal A}({\cal G}(n))\,.
$$
In particular,
 the  differential $d'$ on $E(n,A)$ induces a
 differential, also denoted by $d'$, on $ {\cal A}({\cal G}(n))$:
 $$
 \renewcommand{\arraystretch}{1.6}
 \begin{array}{l} d'(a_1 \otimes a_2 \otimes \ldots \otimes
 a_{l(G)} \cdot e_G) =
 \\
 \mbox{}\hspace{2cm} \sum_{
 \renewcommand{\arraystretch}{0.6}
 \begin{array}{l} {\scriptstyle i<j}\\ {\scriptstyle i,j\in
 S_s}\end{array}
 \renewcommand{\arraystretch}{1}} a_1 \otimes \ldots \otimes a_{l(G)}\cdot
 e_Ge_{ij}
 \\ \mbox{}\hspace{2cm} +         \sum_{
 \renewcommand{\arraystretch}{0.6}\begin{array}{l} {\scriptstyle i<j}\\
 {\scriptstyle i\in S_s\,, j\in S_t}\\{\scriptstyle s<t}\end{array}
 \renewcommand{\arraystretch}{1}}
 (-1)^{\tau_{s,t}} a_1 \otimes \ldots \otimes a_sa_t\otimes \ldots
 \widehat{\otimes a_t}\otimes \ldots \otimes a_{l(G)} \cdot
 e_Ge_{ij}
 \end{array}
 \renewcommand{\arraystretch}{1}
 $$

 For instance, if $n = 3$ and $G$ the discrete graph ($e_G = 1$),
 then
 $$d'(a\otimes b\otimes c) = (ab\otimes c) e_{12} + (a\otimes bc
 )e_{23} + (-1)^{\vert b\vert \cdot \vert c\vert} (ac\otimes b)
 e_{13}\,.$$

It is easy   to prove that  the differential of $A$ induces a
differential $d''$ on the
graded space
${\cal A}({\cal G}(n)) $ such that $(d'+d'')^2=0$. Then Proposition 2.6 of
\cite{BG} extends in:

\vspace{3mm} \noindent {\bf Theorem C.} {\it For any differential graded
commutative algebra
$A$, the isomorphism
$$
\Phi : E(n,A) \to     {\cal A}({\cal G}(n))\,.
$$
is an isomorphism of bicomplexes.}

 \vspace{3mm}
 For any sub-family of graphs $\Gamma \subset {\cal
 G}(n)$ which is stable by adjonction of edges, we       consider the
 subcomplex of ${\cal A}({\cal G}(n))$
$$
{\cal A}\Gamma = \oplus_{G\in \Gamma} A^{\otimes   l(G)}e_G\,.
$$

\vspace{3mm}

We will use the above description of the bicomplex  $E(n,A)$ in order to
replace
it by a smaller
one.

  \vspace{2mm}
\noindent {\bf Proposition 1.} {\sl The ideal  $J$  of $E(n,A)$
 generated by the products $e_{ir}e_{jr}$, $i,j=1, \ldots , n$, is
 a bigraded differential acyclic complex.}

 \vspace{3mm} It follows from Proposition 1 that the canonical  projection

$$
E(n,A) \to \bar E (n,A) := E(n,A)/J
$$         is a quasi-isomorphism.

In fact, as a byproduct of  Proposition 1,   the
isomorphism
$\Phi$  induces an isomorphism of bicomplexes
$$
\bar \Phi : \bar E(n,A) \stackrel{\cong}{\to}    \bar{{\cal A}}({\cal G}(n)):=
\oplus_{G \in
\overline{\cal G}}A^{\otimes l(G)}\, e_G\,,
$$
where $\overline{\cal G}$ denotes the family of graphs of ${\cal
G}(n)$ that do not have two edges with the same target.

 \vspace{2mm}
\noindent {\bf Proof.} Let ${\cal G}_J$ be the set of graphs $G \in {\cal
G}(n)$
such
 that at least two edges terminate at the same vertex. Thus $\Phi$
 induces an isomorphism of differential bigraded vector spaces
 $$
\Phi_J : J \to {\cal A}{\cal G}_J\,.
$$
  We denote by ${\cal
 G}(n,p)$ the subset of ${\cal G}(n)$   consisting of graphs with
       $\leq p$   components, and         by ${\cal G}_J(n, p)$
the intersection         ${\cal
G}_J(n,
 p)=
 {\cal G}_J \cap {\cal G}(n,p)$. We prove by induction on $n$ and
 $p$ that ${\cal A}{\cal G}_J(n, p)$ is acyclic.   We suppose
      that the assertion has been proved   for $q<n$ and any
 $p$, and we consider the disjoint union
$$
{\cal G}_J(n,1) ={\cal H}_0 \cup {\cal   H}_1\cup {\cal H}_2 \,,\,  \left\{
 \begin{array}{l}
 {\cal H}_0 = \{ G \in         {\cal G}_J(n,1) \,\vert \, (1,2)\notin
E(G)\,\}\\
 {\cal H}_1 =\{ G\in         {\cal G}_J(n,1)\,\vert \, (1,2) \in E(G)\,,\,
  l(G{\scriptstyle \backslash} \{(1,2)\}) = 1\,\}\\
 {\cal H}_2 = \{ G\in         {\cal G}_J(n,1)\, \vert \, (1,2) \in E(G) \,,\,
  l(G{\scriptstyle \backslash} \{ (1,2)\}) = 2 \,\}
 \end{array}
 \right.
 $$
   The graded vector space ${\cal A}({\cal H}_1 \cup         {\cal H}_0)$
 is a  sub-complex of ${\cal A}{\cal G}_J(n,1)$. We filter the complex by
putting $e_{ij}$
in degree $i+j$. The first differential
is the internal one. We can therefore  suppose that $d'' = 0$. The next nonzero
differential of the associated spectral sequence, $\delta_3$,
 is the multiplication by $e_{12}$, which corresponds to the addition of the
edge
$(1,2)$. Therefore $\delta_3: \oplus_{G \in         {\cal H}_0}A\, e_G  \to
 \oplus_{G \in       {\cal H}_1}A\, e_G         $ is a linear isomorphism.
       The complex $ {\cal A}{\cal G}_J(n,1)$ is thus quasi-isomorphic to the
quotient
complex
$$
 Q := {\cal A}{\cal G}_J(n,1)/  {\cal A}({\cal H}_0  \cup        {\cal
H}_1)\,.
 $$
 Denote by $Sh(p)$ the set of $(p, n-p-2)$-shuffles of $\{3, 4,
\cdots         ,n\}$,
and for each $(\sigma, \tau)\in Sh(p)$, by ${\cal G}(\sigma
 ,\tau )$  the set of   graphs in ${\cal G}_J$ with two components, the first
one
containing the vertices $1, \sigma_1, \ldots , \sigma_p$
 and the other one the vertices $2, \tau_1, \ldots , \tau_{n-p-2}$.
 Observe now that the map $G \to G{\scriptstyle \backslash} \{ (1,2) \}$
induces
a bijection
 $$
{\cal H}_2 \to \displaystyle\amalg_{p=0}^{n-2}
 \displaystyle\amalg_{(\sigma ,\tau)\in Sh(p)} {\cal G}(\sigma
 ,\tau)\,.
$$

 For $(\sigma ,\tau)\in Sh(p)$ there is clearly a bijection between
 ${\cal G}(\sigma ,\tau )$         and
 $${\cal G}_J \cap \left(\, {\cal G}(p,1) \times {\cal G}(n-p-2,1)\,
 \right)\,.$$         Therefore the quotient complex $Q $ contains the
subcomplex
$$
Q_1:= \oplus_{Sh(p)}\, \left( {\cal A}\, {\cal G}_J(p, 1)\,
 \otimes \,  {\cal A} {\cal G}(n-p-2, 1)\, \right)\,,
$$
that is acyclic since $p<n$.
Moreover the quotient $Q/Q_1$ is isomorphic to
$$
\oplus_{Sh(p)} \left(
{\cal A}\, ( \, {\cal G}(p,1) {\scriptstyle \backslash} {\cal G}_J(p, 1) )
\otimes
{\cal A}{\cal G}_J(n-p-2, 1)\,
\right)\,,
$$
and is acyclic since $n-p-2<n$. Therefore, ${\cal A}{\cal G}_J(n,1)$ is also
acyclic.

 We suppose that ${\cal A}{\cal G}_J(n, p)$ is acyclic for $p<r$,
 and  we consider the quotient complex ${\cal A}{\cal G}_J(n,r) / {\cal
A}{\cal
G}_J(n,r-1)$.

 Denote by $Surj(n,2)$ the set of surjective maps $\{ 1, 2, \cdots
 , n\} \to \{ 1,2\}$. We  associate bijectively to a graph with r
 components in $ {\cal G}_J(n,r) {\scriptstyle \backslash} {\cal
G}_J(n,r-1)$ two graphs $G_1$ and $G_2$, where
 $G_1$ has only one component with $V(G_1) = \pi^{-1}(1)$, and $G_2$
 contains exactly $r-1$ components       with $V(G_2) = \pi^{-1}(2)$.
 This produces  an isomorphism of complexes
 $$
 \renewcommand{\arraystretch}{1.6}
 \begin{array}{l}   {\cal A}{\cal G}_J(n,r) / {\cal A}{\cal G}_J(n,r-1))=   \\
  \oplus_{\pi \in Surj(n,2)} (       {\cal A}\,
 {\cal G}_J (\pi^{-1}(1), 1) \otimes         {\cal A}{\cal G}(\pi^{-1}(2),
r-1))
\\
 \mbox{}\hspace{2cm}         \oplus {\cal A}\, ({\cal G}(\pi^{-1}(1), 1)

 {\scriptstyle \backslash}{\cal G}_J (\pi^{-1}(1), 1))
 \otimes {\cal A}{\cal G_J}(\pi^{-1}(2), r-1)))  \,.
  \end{array}
\renewcommand{\arraystretch}{1}
 $$
 The induction hypothesis shows that this quotient complex is
 acyclic.       \hfill $\square$

 \section{Proof of Theorem 1}

\vspace{3mm} Let $H = H^*(M;\bk)$. To avoid confusion we denote by $T(n,H)$
 the $E_1$ term of the Cohen-Taylor spectral sequence,
$$
T(n,H) = H^{\otimes n}\otimes \land (x_{ij})/I\,.$$
We consider   the sub-vector space, $R$,  of $\land x_{ij}$ generated by
the monomials
$$x_{i_1,j_1}\cdots x_{i_r,j_r} \hspace{1cm}\mbox{with }
\hspace{3mm}\left\{ \begin{array}{l}
0\leq r\leq n\\
i_s<j_s\,, \mbox{ for $s=1, \ldots ,r$}\\
j_1 <\cdots <j_r\,.
\end{array}
\right.
$$
 The differential $H^{\otimes n}$-module
$H^{\otimes n}\otimes R$ is a direct summand of $H^{\otimes n}\otimes \land
(x_{ij})$. Let $L$ be      the sub $H^{\otimes
n}$-module of $H^{\otimes n} \otimes R$  generated by the elements
$(e_i-e_j)\otimes x_{ij}$, $i<j$. The
quotient  map
$(H^{\otimes n}\otimes R)/L \to T(n,H)$ is a bigraded isomorphism of
$H^{\otimes n}$-modules.

 Since  $M$ is a connected closed $m$
dimensional manifold  one can choose
$\omega
\in H^m$ such that:

  a) $H^m=\bk \omega$

  b) the relations $ab = <a;b> \omega\,, a \in H^p\,, b\in H^q $ define  non
degenerated
bilinear forms $ H^p \otimes H^{n-p} \to \bk \,, \quad a\otimes b \mapsto
<a;b>
$.

This pairing extends to a pairing
$$
(\ast) \qquad H^{\otimes n} \otimes H^{\otimes n} \to \bk \,, \,
< a_1\otimes...\otimes a_n; b_1\otimes ... \otimes b_n=
\epsilon(\sigma) < a_1;b_{1} > ... < a_n ; b_{_n} >
$$
 where $\epsilon (\sigma)$ is the graded signature.

In the same fashion, the pairing
$$
\left( \oplus _{i<j} x_{ij}\bk  \right) \otimes \left( \oplus _{i<j} s^{1-m}
e_{ij}\bk \right)\to
\bk
\,, \, <x_{ij}, s^{1-m}e_{kl}> = \delta _i^l \delta _j^k
$$
 extends to a pairing
$$
(\ast \ast) \qquad  \land (x_{ij}) \otimes \land(s^{1-m}e_{kl}) \to \bk
\,,\quad
1\leq i<j\leq
n\,,
\,1\leq k<l\leq n\,.
$$
The pairings $(\ast)$ and $(\ast \ast)$ glue together  in the pairing
 $$
 \qquad \left( H^{\otimes n} \otimes \land (x_{ij}) \right) \otimes \left(
H^{\otimes
n}\otimes
\land(s^{1-m}e_{kl})\right) \to \bk \,,\quad 1\leq i<j\leq
n\,,
\,1\leq k<l\leq n \,.
 $$
which in turn defines the $H$-linear isomorphism
$$
  H^{\otimes n} \otimes \land (x_{ij})   \to \left(H^{\otimes
n}\otimes
\land(s^{1-m}e_{kl})\right)^\vee\,, \qquad z \mapsto <z;-> \,.
$$

 This isomorphism produces a  $H^{\otimes n}$-linear isomorphism
$$
\theta  : (H^{\otimes n}\otimes R)/L \rightarrow \bar E(n,H)
$$
which restricts, for each $p\geq 0$, to  the isomorphisms
$$
\theta^{p,q} :   \left( (H^{\otimes n}\otimes R)/L\right)
^{p,q}
\stackrel{\cong}{\to}
\left(\bar  E(n,H)\right) ^{p,q-pm}\,.$$

Recall now that the differential $d_1$ on $T(n,H)$ is defined by
$d_1(x_{ij}) = p_{i,j}^*(\delta^M)$, where $\delta^M$
denotes the diagonal class. Since $M$ is compact $$\delta^M = \sum_t
(-1)^{\vert e_t'\vert} e_t\otimes e_t'$$
where $\{e_t\}$
denotes a linear basis of $H^*(M;\bk)$  and $\{e_t'\}$ the Poincar\'e dual
basis.

\vspace{3mm}
Theorem 1 is then a direct consequence of the commutativity of the
following diagram.
$$
\begin{array}{cccc}
(T(n,H)) ^{p,\ast}& \stackrel{\theta^p   }{ \longrightarrow}& \left(
\bar E(n,H)^{p,\ast}
\right) ^\vee
\\ {\scriptstyle d_1} \downarrow && \uparrow{\scriptstyle  (d_1)^\vee }\\
(T(n,H)) ^{p+1,\ast}& \stackrel{ \theta^{p+1}  }{ \longrightarrow} &
\left( \bar E(n,H)^{p+1,\ast}\right) ^\vee\,,
\end{array}
$$
where $ {\bar E}^\vee:=Hom
(\bar E,\bk)$ denotes
the graded dual of the graded vector space $\bar E$.
\hfill $\square$

 \section{The bicomplex $C(n,A)$}

Let $A$ be a commutative differential graded algebra. Recall from section 1
the isomorphism  of bicomplexes
$$
\bar \Phi : \bar E(n,A) \stackrel{\cong}{\to}    \bar{{\cal A}}({\cal G}(n)):=
\oplus_{G \in
\overline{\cal G}}A^{\otimes l(G)}\, e_G\,,
$$
where $\overline{\cal G}$ denote the family of graphs of ${\cal
G}(n)$ that do not have two edges with the same target.
 Observe that each
 component of $G \in \overline{\cal G}$ is a tree and that the
 number of edges is exactly $n-l(G)$.

 We denote by ${\cal H}$ the
 subset of $\overline{\cal G}$ consisting of graphs such that
 $S_{1} = \{ 1\}$.   We suppose that $A = \bk \oplus A^+$, with $A^+
 = \oplus_{i>0}A^i$, and we consider the new bicomplex
 $$
 C(n,A) = \left(\oplus_{G \in {\cal H}} A \otimes (A^+)^{\otimes
 l(G)-1)}e_G, d\right)
$$
 whose differential $d = d' + d"$ is defined as follows: $d'':
C(n,A)^{\ast, k }
\to
C(n,A)^{\ast , k+1} $ is
 the internal differential coming from the differential on $A$, and
 $$
d': C(n,A)^{k,\ast } \to C(n,A)^{k+1,\ast } \,, \quad  d'(a_1\otimes \ldots
\otimes a_{l}e_G)
=
\sum_{1< i<j\leq n}
 \alpha_{ij}\, e_G\, e_{ij}\,,
$$ with the conditions
 \begin{enumerate}
 \item[$\bullet$] the term $e_G\,e_{ij}$ is zero if $G\cup \{(i,j)\}\not\in
{\cal
H} $;
 \item[$\bullet$] In the other case, denote by $S_s$ and $S_t$, $1<s<t$, the
 components of $i$ and $j$, then $G' = G\cup \{i,j\}$ has $l-1$
 components. The element $\alpha_{ij}$ is the element of $A
\otimes (A^+)^{\otimes l-2}$ defined by
$$
\begin{array}{l}
 \alpha_{ij}=(-1)^{\vert a_t\vert (\vert a_{s+1} \vert + \ldots \vert
 a_{t-1}\vert)} \, a_1 \otimes a_2 \ldots \otimes a_{s-1}\otimes
 a_sa_t \otimes a_{s+1} \otimes \ldots \widehat{\otimes a_t}\otimes
 \ldots \otimes a_l
 \\
 \mbox{}\hspace{1cm} - \varepsilon \cdot a_1 a_s \otimes a_2 \ldots
 a_{s-1}\otimes a_t \otimes a_{s+1} \otimes \ldots
 \widehat{\otimes a_t}\otimes \ldots \otimes a_l \\
 \mbox{}\hspace{1cm} - (-1)^{\vert a_t\vert(\vert a_2\vert + \ldots
 \vert a_{t-1}\vert ) }\, a_1a_t\otimes a_2\ldots \otimes
 a_{s-1}\otimes a_s \otimes a_{s+1} \otimes \ldots \widehat{\otimes
 a_t}\otimes \ldots \otimes a_l\,,
 \end{array}
 $$
 with $\varepsilon = (-1)^{\vert a_s\vert(\vert a_2\vert         + \ldots
 \vert a_{s-1}\vert )         + \vert a_t\vert (\vert a_{s+1} \vert +
 \ldots \vert a_{t-1}\vert)}$.
 \end{enumerate}

 \noindent For instance, for $n = 3$ and $G$       the discrete graph, we have
 $$d'(a_1\otimes a_2\otimes a_3) = \left[ a_1 \otimes a_2a_3 -
 a_1a_2\otimes a_3 - (-1)^{\vert a_2\vert \cdot \vert a_3\vert }
 a_1a_3\otimes a_2\,\right] \, e_{23}\,. $$

 \vspace{3mm}

\noindent{\bf Proposition 3.} {\sl
 There exists a quasi-isomorphism         of bicomplexes
$$     \bar \varphi :     C(n,A)
 \stackrel{\simeq}{\to} \bar E (n,A)
 $$ which is natural in $A$.}

 \vspace{3mm}\noindent {\bf Proof.} Let us consider for each $l\geq
 2$ the linear map
 $$\gamma_l : A^{\otimes l} \to A^{\otimes l}$$
 defined by
 $$\gamma_l (a_1, \ldots , a_{l}) = \sum         (-1)^{k} \varepsilon (i_1,
 \ldots , i_k) \, a_1 a_{i_1}\ldots a_{i_k} \otimes         {a_2'} \otimes
 \cdots \otimes         {a_{l}'}\,,$$ where $ {a_j'}= a_j$ if $j\not\in \{
 i_1, \ldots , i_k\} $ and equals $1$ otherwise. The sum is taken
 over all the ordered (possibly empty) subsets
  $i_1<i_2<\ldots < i_k  $
 of $\{ 2, \ldots , l\}$, and $\varepsilon (i_1, \ldots , i_k)$ is
 the graded signature  of the permutation
  $
a_1 a_2 \ldots a_l \mapsto
 a_1a_{i_1}\ldots a_{i_k} {a_2}'\ldots {a_l}'\,.
 $         We define
 $$\varphi : C(n,A) \to  {\cal A}{\cal G}(n)\,, \hspace{5mm}
\varphi (a_1\otimes \ldots \otimes a_{l}\, e_G) =\gamma_l
(a_1,
\ldots , a_{l})\, e_G \,,$$
which is an injective linear map that commutes
 with the differentials.

 Now observe that the composite $\bar \varphi : C(n,A)
 \stackrel{\Phi^{-1}\circ \varphi}{\to} E(n,A) \to \bar E(n,A)$ remains
injective
 and that for any $r\geq 2$
 $$
\mbox{Im} \bar \varphi \cdot e_{1r} = 0\,.
$$
 Let denote by $(\Gamma, \bar d)$ the cokernel of $\bar\varphi$:
 $$(\Gamma, \bar d) = (\, \displaystyle\oplus_{G \in \overline{\cal G}} \,
 C(G)e_G,\bar d\, )\,,$$
 with $$ C(G) = \left\{\begin{array}{l} A^{\otimes
 l}\,,\hspace{3mm} \mbox{
 if $G$ contains an edge $(1,r)$ for some $r$}\\         \oplus_{s=1}^{l-1}
A^{\otimes s} \otimes \bk \otimes
 A^{\otimes l-s-1}\,, \hspace{3mm}\mbox{otherwise} \end{array}
 \right.$$

 The differential $\bar d$  induced by $d$   is the sum $ \bar
 d = d" + \bar d_0 + \bar d_1$, where $d"$ is the internal
 differential, $\bar d_0$ consists into the multiplication by
 $\sum_s e_{1s}$ and $\bar d_1$ comes from the multiplication by
 $\sum_{2\leq r<s} e_{r,s}$.

 We prove that $H^*(\Gamma, \bar d_0) = 0$, which implies the
 result by an elementary spectral sequence argument. Let ${\cal H}$
 defined as before, and         consider         the partition
 $$\overline{\cal G} = {\cal H} \cup {\cal G}_2 \cup {\cal G}_3 \cup
 \ldots \cup {\cal G}_n\,,$$ where         ${\cal G}_r, r\geq 2$, denotes
 the set of graphs $G$ such that $r$ is the smallest vertex $\neq
 1$ in the connected component of $1$. We       denote by ${\cal G}_r'$
 the set of graphs $G$ in $\overline{\cal G}$ such that         the
 components of $1$ and $r$ are different and do not contain any of
 the vertices $2, \ldots , r-1$. This implies a direct sum
 decomposition of   $\Gamma$:
 $$\Gamma = \oplus_{r=2}^n (\Gamma_r \oplus \Gamma_r')\,,$$
 where $\Gamma_r = \oplus_{G \in {\cal G}_r} \, D(G)\, e_G$ with
 $D(G) = D_1(G) \otimes \ldots \otimes D_l(G)$ and $$ D_i(G) =
 \left\{
 \begin{array}{ll} A^+ & \mbox{ if         $S_i$ contains one of the
 vertex $2, \ldots , r-1$}\\A         & \mbox{otherwise,}
 \end{array}
 \right.
 $$
 and where, $\Gamma_r' = \oplus_{G \in {\cal G}_r'} D'(G)e_G$, with
 $D'(G) = D_1'(G) \otimes \ldots \otimes D_l'(G)$, and
 $$D_i'(G) = \left\{
 \begin{array}{ll}
 A^+ & \mbox{ if         $S_i$ contains one of the vertex $2, \ldots
 , r-1$}\\
 \bk & \mbox{if $S_i$ contains the vertex $r$}\\         A &
\mbox{otherwise.}
\end{array}         \right.         $$

 To establish the decomposition $\Gamma = \oplus_{r=2}^n (\Gamma_r
 \oplus \Gamma_r')\,$, let $a = a_1 \otimes \ldots \otimes a_{l}\,
 e_G$ be an element in         $\Gamma$. When $G \in {\cal H}$, we denote
 by $r$ the smallest       vertex contained in a component $S_i$ such
 that $a_i\in k$. In this case $a \in \Gamma_r'$. Suppose now that
 $G         \in {\cal G}_r$. If for each $i$, $2 \leq i <r$, the vertex
 $i$ belongs to a component $C_s$ with $a_s\in A^+$, then $a\in
 \Gamma_r$. In the other case, denote by $t$ the smallest vertex
  in a component $C_s$ with $a_s         \in \bk$, then $a \in
 \Gamma_t'$.

   The vector space $\Gamma_2 \oplus \Gamma_2'$ is a
subcomplex of $\Gamma$,         and the multiplication by $e_{12} $ induces
an
isomorphism $\Gamma_2' \stackrel{\cong}{\to} \Gamma_2$; in the quotient
 $\Gamma /         (\Gamma_2 \oplus \Gamma_2')$, the vector space $\Gamma_3
\oplus         \Gamma_3'$ is a subcomplex and the multiplication by
$e_{13}$ is
an isomorphism         $\Gamma_3' \stackrel{\cong}{\to} \Gamma_3$. More
generally,         $\Gamma_r \oplus \Gamma_r'$ is a subcomplex in the
quotient $
\Gamma /         ((\oplus_{s<r}\Gamma_s )\oplus (\oplus_{s<r}\Gamma_s'))$, the
vector space
 $\Gamma_r \oplus         \Gamma_r'$ is a subcomplex and the multiplication by
$e_{1r}$ is         an isomorphism         $\Gamma_r' \stackrel{\cong}{\to}
\Gamma_r$. This shows that         $(\Gamma, \bar d_0)$ is acyclic.
\hfill $\square$

 \vspace{3mm}     Theorem 2 of the Introduction follows   directly from the
next Proposition

 \vspace{3mm}\noindent {\bf Proposition 4.} {\sl  Let $M$ be a simply
connected compact oriented manifold. Then $E_{n-1} =E_{\infty}$ in
 the Bendersky-Gitler spectral sequence for $H^*(F(M;n);\bk)$.}

 \vspace{3mm}\noindent {\bf Proof.}         The $E_2$ term   of the
 Bendersky-Gitler spectral sequence is isomorphic as a bigraded
 vector space to the $E_2$ term  of the spectral sequence associated
 to the bicomplex $C(n,A)$. Since   $C(n,A)^{\geq n-1,*} = 0$, we get
$E_2^{p,q} =
0$ for $p>n-2$.
\hfill $\square$

 \vspace{3mm}

  If $H$ is a  differential graded commutative algebra then $C(3,H)$ is the
complex ($d"=0$)
 $$
(H \otimes H^+ \otimes H^+) \stackrel{d_1}{\longrightarrow}
 (H\otimes H^+)\cdot e
$$
where $e$ is a variable of degree 1, and
 $$ d_1(a\otimes b\otimes c) = a\otimes  bc - ab \otimes c - (-1)^{\vert
c\vert
\cdot \vert b\vert}
ac
 \otimes b\,.$$
When  $H = H^*(M;\bk)$, the kernel of $d_1$, $E_2^{0,*}$, is the image
$I_M$ of the
map
 $$H^*(M^3;\bk )\to H^*(F(M,3);\bk)$$
 induced by the inclusion $F(M,3) \hookrightarrow M^3$. The
 cokernel of $d_1$ is   the first Hochschild homology group of
 the unital algebra $H $ with coefficients in the  bimodule $H$,
  $H\!H_1(H;H)$, and by
 (\cite{Lo}, 1.1.10)   is    isomorphic to the $H$-module of         K\"ahler
 differentials $\Omega^1_{H/\bk}$ on $H$. It follows from the
 Lefchetz duality $H^p(M^n,D^nM) \cong H_{nm-p}(F(M,n);\bk)$ that

 \vspace{3mm}\noindent {\bf Proposition 5.} {\sl For a
 $m$-dimensional oriented closed manifold $M$ there is a short
 exact sequence
 $$0 \to (I_M)^* \to H^*(F(M,3);\bk) \to       (\mbox{Hom}
 (\Omega^1_{H/\bk}, \bk))^{*+1-3m} \to 0\,.$$         }

 \section{  $F(M,4)$ and Massey products}

 For a commutative differential graded algebra $(A,d)$, the double
 complex $C(4,A)$ has the form
 $$0\to A\otimes (A^+)^{\otimes 3} \stackrel{d'}{\longrightarrow}
C(4,A)^{1,\ast}
 \stackrel{d'}{\longrightarrow} C(4,A)^{2,\ast}\to 0\,,$$ with
$$
 \renewcommand{\arraystretch}{1.6}
 \begin{array}{l}
 C(4,A)^{1,\ast} = (A\otimes A^+\otimes A^+)e_{23} \oplus (A\otimes
 A^+\otimes A^+)e_{24}
 \oplus (A\otimes A^+\otimes A^+)e_{34} \,,\\
       C(4,A)^{2,\ast} = (A\otimes A^+) e_{23}e_{24} \oplus (A\otimes
A^+)e_{23}e_{34}\,
\end{array}         \renewcommand{\arraystretch}{1}
$$
$$
\renewcommand{\arraystretch}{1.5}
\begin{array}{l}         d'(a\otimes b\otimes c\otimes
d) = (\, a\otimes bc\otimes d - ab\otimes c
 \otimes d -         (-1)^{\vert b\vert\cdot \vert c\vert} ac\otimes b \otimes
d\, ) \, e_{23}
 \\   \mbox{}\hspace{5mm} + (\, a\otimes b\otimes cd - (-1)^{\vert b\vert\cdot
 \vert c\vert} ac\otimes b\otimes d         - (-1)^{\vert d\vert (\vert b\vert
+ \vert c\vert )} ad\otimes b\otimes
 c\, )\,         e_{34}\\   \mbox{}\hspace{5mm}
 + (\, (-1)^{\vert d\vert \cdot \vert c\vert
}
 a\otimes bd\otimes c - (-1)^{\vert d\vert \cdot \vert c\vert } ab\otimes
d\otimes
c         -(-1)^{\vert d\vert (\vert b\vert +\vert c\vert )} ad\otimes
b\otimes
c\, )\,
 e_{24}\,,
 \\\mbox{}
 \\         d'((a\otimes b\otimes c) \, e_{23}) = \Delta (a,b,c) \,
(e_{23}e_{34}+
 e_{23}e_{24})\\ d'((a\otimes b\otimes c)\, e_{24}) = - \Delta (a,b,c) \,
e_{23}e_{24}\\         d'((a\otimes b\otimes c)\, e_{34}) = -\Delta (a,b,c) \,
e_{23}e_{34}\,.
 \end{array}
 \renewcommand{\arraystretch}{1}$$
 where $\Delta(a,b,c) = (a\otimes bc - ab\otimes c-(-1)^{\vert
 b\vert \cdot \vert c\vert } ac\otimes b)$

 \vspace{3mm}  From Theorem B we deduce

 \vspace{3mm}\noindent {\bf Proposition 6.}         {\sl Let $\bk$ be a
field of characteristic zero, the $E_2$
  term of the Bendersky-Gitler spectral sequence for
 $H^*(F(M,4);\bk)$ satisfies
$$E_2^{2,*} = \Omega^1_{H/\bk}\,
 e_{23}e_{24} \, \otimes \, \Omega^1_{H/\bk}\ e_{23}e_{34}\,,$$
with $H = H^*(M;\bk)$.
 }

 \vspace{3mm} In order to detect nonzero element in $E_2^{2,*}$, we
 now describe a special quotient of $E_2^{2,*}$.       Denote $H =
 H^*(M;\bk)$ and $Q(H) = H^+/(H^+\cdot H^+)$, the space of indecomposable
elements, and let
 $$\psi : Q(H)\otimes Q(H) \to Q(H)\otimes Q(H)$$ be the linear map
 defined by $\psi (a\otimes b) = a\otimes b - (-1)^{\vert a\vert
 \cdot \vert b\vert} b\otimes a$. The image of $\psi$ is the sub-vector
space of $\Sigma_2$-invariants elements in $Q(H)\otimes
 Q(H)$, where the action of the generator $\tau$ of $\Sigma_2$ is
 defined by
 $$\tau (a\otimes b) = -(-1)^{\vert a\vert \cdot \vert b\vert}
 b\otimes a\,.$$ We now consider the composition
 $$
 \begin{array}{ll}
 H \otimes H^+ \stackrel{\cong}{\to} (\bk \otimes H^+ ) \oplus
 (H^+\otimes H^+) & \to (\bk \otimes Q(H) ) \oplus (Q(H)\otimes
 Q(H))
 \\ & \stackrel{\psi}{\to} (\bk \otimes Q(H) ) \oplus (Q(H)\otimes
 Q(H))^{\Sigma_2}\,. \end{array}
 $$
 By definition of $d'$, this morphism induces a surjective map from
 $ H^{2,\ast}(C(4,H), d_1)$ onto $$ \left[(\bk \otimes Q(H) ) \oplus
 (Q(H)\otimes Q(H))^{\Sigma_2}\right]\,\, e_{23}e_{34} \, \oplus\,
 \left[(\bk \otimes Q(H) ) \oplus (Q(H)\otimes
 Q(H))^{\Sigma_2}\right]\,\,e_{23}e_{24}\,.$$

 \vspace{3mm} Recall now the definition of a triple matrix Massey
 product (\cite{M}). Consider   three matrices $L, B, C$
 $$L = (a_1, \ldots , a_r)\,,\hspace{1cm} B =
 \left(\begin{array}{ccc} b_{11} &\ldots & b_{1s}\\ \ldots &\ldots
 &\ldots\\
 b_{r1}&\ldots & b_{rs}
 \end{array}
 \right)\,,\hspace{1cm} C = \left(\begin{array}{c} c_1\\ \ldots \\
 c_s\end{array}\right)\,,$$ whose entries are       cocycles in $A$       such
 that
 $$L\cdot B = (d(x_1), \ldots , d(x_s))\,\hspace{5mm}\mbox{ and
 }\hspace{5mm} B\cdot C =
 \left(\begin{array}{c} d(y_1)\\\vdots \\
 d(y_r)\end{array}\right)\,.$$         The triple matrix Massey product
 $\langle L,B,C\rangle$ is defined, up to some indeterminacy, as the
 class of the cocycle $\sum_i x_ic_i - \sum_j (-1)^{\vert a_j\vert}
 a_jy_j$. Nonetheless the residual class of         $\langle L,B,C\rangle$ in
the quotient $Q(H^*(A))$ is
 uniquely defined.

 \vspace{3mm}\noindent {\bf Theorem 4.} {\sl  Let $\bk$ be   of
characteristic zero, $M$ be   simply
 connected  and $N = M\# ( S^2 \times S^{m-2}) $ be the connected sum of $M$
with
 a product  of two  spheres.  If
  there is  a triple matrix Massey
 product in $H^*(M;\bk)$ whose residual class in $Q(H^*(M;\bk))$ is
nonzero, then the Cohen-Taylor spectral
 sequence for   $H^*(F(N,4);\bk)$ does not collapse at the $E_2$-term.  }

 \vspace{3mm}\noindent {\bf Proof.}  Let us take a Sullivan model
$(A,d)$ for $M$ over the field $\bk$ that satisfies the  following
properties: $A$ is finite
dimensional, $A^0 =
\bk$, $A^1 = 0$, $A^m = \bk\omega$, and $A^{>m} = 0$.
 We choose for the product $S^2 \times S^{m-2}$ the model $(\land
 (x,y)/ (x^2, y^2), d= 0)$ with $x$ in degree 2 and $y$ in degree
 $m-2$. Therefore, (\cite{FHT}), a model for the connected sum $N =
 (S^2\times S^{m-2})\# M$ is given by the differential graded
 algebra
 $$(\left[ A\times_{\bk} \land(x,y)/(x^2,y^2)\right] / \omega
 - xy, d)\,.$$

 Suppose that $A $ admits a triple matrix Massey
 product $\langle L,B,C\rangle$ that is non zero in
 $H^+(M)/(H^+(M)\cdot H^+(M))$, then with the above notation for $L$, $B$
and $C$, we consider the element
 $$
u = \sum_{ij} x\otimes a_i\otimes b_{ij}\otimes c_j -
 \sum_{ij}(-1)^{\vert c_j\vert \cdot \vert b_{ij}\vert + \vert
 c_j\vert \cdot \vert a_i\vert + \vert b_{ij}\vert \cdot \vert
 a_i\vert} x\otimes c_j\otimes b_{ij}\otimes a_i\,.
$$
Then
 $d''u =0$  and $d_1([u])=0$ in the  first spectral sequence of the
bicomplex $C(4,A)$.
 More precisely,
 $$\renewcommand{\arraystretch}{1.5}
 \begin{array}{rl}
 d'(u) = & \left(\sum_jx\otimes d(x_j)\otimes c_j - (-1)^{         \vert
 c_j\vert \cdot \vert a_i\vert + \vert b_{ij}\vert \cdot \vert
 a_i\vert}\sum_i x\otimes d(y_i)\otimes a_i\right)\, e_{23} \\ +&
 \left(\sum_i (-1)^{\vert a_i\vert } x\otimes a_i\otimes d(y_i) -
 (-1)^{\vert c_j\vert +\vert c_j\vert \cdot \vert b_{ij}\vert +
 \vert c_j\vert \cdot \vert a_i\vert } \sum_j x\otimes c_j\otimes
 d(x_j)\right)\, e_{34}
 \end{array}
 \renewcommand{\arraystretch}{1}
 $$
 Therefore
 $$d_2([u]) = \left[x\otimes \langle L,B,C\rangle\right]\,\,
 e_{23}e_{24} \, + 2 \left[x\otimes \langle
 L,B,C\rangle\right] \,\, e_{23}e_{34}\,.$$ The element
 $d_2([u])$ is nonzero in $E_2$ because its image is nonzero
 in
 $$\left[\left(\, Q(H) \otimes Q(H)\,
 \right)^{\Sigma_2}\right]\,e_{23}e_{24}\,.$$ \hfill $\square$

 \vspace{3mm}
Consider now the case where $a,b,c$ and $d$ are elements in $H^*(M;\bk)$, such
that   $ab
= ac= ad= bc=bd=cd= 0$. In the first spectral sequence of the bicomplex
$C(4,H)$
 we obtain the         formula
 $$
 \renewcommand{\arraystretch}{1.6}
 (*) \hspace{3mm}
 \begin{array}{l}
d_2([a\otimes b\otimes c\otimes d])         = \left[ \right.
 (-1)^{\vert a\vert} a\otimes \langle b,c,d\rangle       + \langle
 a,b,c\rangle \otimes d \\         \mbox{}\hspace{5mm}         + (-1)^{\vert
 c\vert\cdot \vert a\vert \cdot \vert b\vert} \langle b,a,d\rangle
 \otimes c \left. -(-1)^{\vert b\vert \cdot \vert c\vert + \vert
 b\vert \cdot \vert d\vert + \vert c\vert\cdot \vert d\vert}
 \langle a,d,c\rangle \otimes b\, \right]\,\,
 e_{23}e_{34}\\   +       \left[ \right.       (-1)^{\vert a \vert +
 \vert b\vert\cdot \vert c\vert} a\otimes \langle c,b,d\rangle         +
 (-1)^{\vert b\vert \cdot \vert c\vert } \langle a,c,b\rangle
 \otimes d\\
 \mbox{}\hspace{5mm}+ (-1)^{\vert b\vert \cdot \vert c\vert + \vert
 b\vert \cdot \vert d \vert + \vert a \vert \cdot \vert c \vert}
 \langle c,a,d\rangle \otimes b         \left. -(-1){ \vert b\vert \cdot
 \vert d\vert + \vert d\vert \cdot \vert c\vert} \langle
 a,d,b\rangle \otimes c\, \right]\,\, e_{23}e_{24}
 \end{array}
 \renewcommand{\arraystretch}{1}
 $$

  \vspace{3mm} We deduce:

 \vspace{3mm}\noindent {\bf Theorem 3.} {\sl  Let $\bk$ be   of
characteristic zero and $M$ be   simply
 connected. Suppose there exist indecomposable elements  $a,b,c$
and $d\in
H^*(M;\bk)$, such that
 \begin{enumerate}
 \item[$\bullet$] $ab = ac= ad= bc=bd=cd= 0$,
 \item[$\bullet$]  the triple Massey product $\langle b,c,d\rangle $ is a
 nonzero indecomposable element in cohomology, and
 \item[$\bullet$] $a$ is not a linear
 combination of the elements $b,c,d$
and $\langle b,c,d\rangle $,
 \end{enumerate}  then
the Cohen-Taylor spectral
 sequence for   $H^*(F(M,4);\bk)$ does not collapse at the $E_2$-term. }

 \vspace{3mm}\noindent {\bf Proof.} By formula (*), the image of
 $d_2([a\otimes b\otimes c\otimes d])$ in $(Q(H) \otimes
 Q(H))^{\Sigma_2}$ is nonzero because the element $a$ is not a
 linear combination of the elements         $b,c,d$ and $\langle
 b,c,d\rangle$.

 \hfill $\square$

       \vspace{3mm}\noindent {\bf Example.} Let $M$ be the sphere tangent
bundle
to
the manifold
 $S^2\times S^2$. A Sullivan minimal model for $M$ is given by
 $(\land (x,y,u,v,t),d)$ with $\vert x\vert =\vert y\vert = 2$,
 $\vert u\vert =\vert v\vert =\vert t\vert = 3$, $d(u)         = x^2$,
 $d(v) = y^2$ and $d(t) = xy$ (\cite{FHT}). A basis of the
 cohomology is given by $1, [x], [y], [tx-uy], [ty-vx], [txy -
 uy^2]$. The cohomology classes in degree 5 are usual triple Massey
 products represented by undecomposable elements in cohomology.
 Formula (*) gives
 $$
d_2([x]\otimes [x]\otimes [y]\otimes [y])
 = ([x]\otimes [ty-vx] - [ty-vx]\otimes [x] -2 [tx-uy] \otimes
 [y])e_{23}e_{24} + \mu e_{23}e_{34}\,.
$$
 The component of $e_{23}e_{24} $ is not symmetric, hence not zero in
$E_2^{2,*}$.
The Bendersky-Gitler spectral sequence for
 $H_{mn-\ast}(F(M,4);{\mathbb Q})$
 does not collapse at the $E_2$ term .

 \vspace{1cm}

    D\'epartement de Math\'ematique Universit\'e Catholique de Louvain
2, Chemin du
Cyclotron   1348  Louvain-La-Neuve, Belgium

 \vspace{3mm} Facult\'e des Sciences Universit\'e d'Angers 2, Boulevard
Lavoisier 49045 Angers,
France


\begin{thebibliography}{xx}


\bibitem{An} D. Anick,  Hopf algebra up to homotopy, {\sl J. Amer.
Math. Soc.}  {\bf 2} (1989) 417-453.

\bibitem{BG} M. Bendersky and S. Gitler, The cohomology
of certain function spaces, {\sl Trans. Amer. Math. Soc.} {\bf 326} (1991),
423-440.

\bibitem{Bo} C.-F. Bodigheimer, F. Cohen and L. Taylor, On the
homology
of configuration spaces, {\sl Topology}         {\bf 28}   (1989), 111-123.

\bibitem{C} F.
Cohen and L. Taylor, Computations of Gelfand-Fuks   cohomology, the cohomology
of function spaces,
and the   cohomology of configuration spaces, Lecture Notes in Math. {\bf
657}, Springer-Verlag, 1978, 106-143.

\bibitem{Del} P. Deligne, P. Griffiths, J. Morgan, and D. Sullivan,
       Real homotopy types of K\"alher manifolds, {\sl Invent. Math.}
       {\bf 29} (1975) 245-274.

\bibitem{El} M. Elhaouari, {  p-formalit\'e des espaces}, {\sl J. of Pure and
       Applied Algebra} {\bf 78} (1992) 27-47.

 \bibitem{FT} Y. Felix and J.-C. Thomas, {\sl Rational homotopy
type of
 configuration spaces}, Topology and its applications {\bf 102}
 (2000), 139-149.

\bibitem{FHT}
 Y. Felix, S. Halperin and J.-C. Thomas, {\sl Rational Homotopy
 Theory}, Graduate Text in Mathematics, Springer Verlag, 2000.

\bibitem{K} I. Kriz, On the rational homotopy type of
 configuration spaces, {\sl Ann. Math.}         {\bf 139} (1994),   227-237.

\bibitem{Lo} J.L.
Loday, {\sl Cyclic homology}, Grundlehren der   Mathematischen Wissenschaften
{\bf 301},
Springer-Verlag, 1992.

 \bibitem{LoM} P. L\"offler
and J. Milgram, The structure of   deleted symmetric products,         in {\sl
Braids}, Contemporary   Mathematics {\bf 78}, American Mathematical Society
(1988),   415-424.

\bibitem{M} J.P. May, Matrix Massey products, {\sl J. of   Algebra} {\bf
12} (1969), 533-568.

\bibitem{Su} D. Sullivan, Infinitesimal computations in   topology, {\sl Publ.
IHES} {\bf         47 } (1977), 269-331.

\bibitem{T} B. Totaro, Configuration spaces of algebraic
 varieties, {\sl Topology} {\bf 35} (1996), 1057-1067.

 \end{thebibliography}
\end{document}